\documentclass[11pt]{article}
\usepackage{latexsym, amscd, amsmath, graphics, rotating, array}
\usepackage{amssymb,amsmath,multibox}

\makeindex
\title{Counting subgroups in a family of nilpotent semi-direct products}
\author{Christopher Voll\thanks{Mathematical Institute, Oxford OX1 3LB, England.}}

\newenvironment{acknowledgements}{\bigskip\noindent\textsl{Acknowledgements.}\rm}
\newenvironment{proof}{\bigskip\noindent\textbf{Proof.}\rm}{\hfill$\Box$ \\ \hspace{0.1cm}}

\newtheorem{lemma}{Lemma}
\newtheorem{theorem}{Theorem}

\newtheorem{corollary}{Corollary}
\newtheorem{proposition}{Proposition}

\def \N {\ensuremath{\mathbb{N}}}

\def \isom {\ensuremath{\stackrel{\sim}{=}}}

\def \Qp {\mathbb{Q}_p}
\def \Z {\mathbb{Z}}

\def \Zp  {\mathbb{Z}_p}

\def \Fp  {\mathbb{F}_p}

\def \LieL {\ensuremath{\mathfrak{L}}}

\def \F23 {\ensuremath{F_{2,3}}}

\def \wnjxy {W_{n,J}({\bf X},{\bf Y})}

\def \bfX {{\bf X}}
\def \bfY {{\bf Y}}
\def \LieLp {\LieL_p}
\def \ul {\underline}

\begin{document}
\maketitle

\begin{abstract}
\noindent In this paper we compute the subgroup zeta functions of
nilpotent semi-direct products of groups of the form $$G_n := \langle
x_1,\dots,x_{n},y_1,\dots,y_{n-1}|\;[x_i,x_n]=y_i,\;1\leq i \leq
n-1,\text{ all other $[,]$ trivial}\rangle$$ and deduce local functional equations.
\end{abstract}

\section{Introduction}
 For $n\geq 2$ define the group $$G_n := \langle
x_1,\dots,x_{n},y_1,\dots,y_{n-1}|\;[x_i,x_n]=y_i,\;1\leq i \leq
n-1,\text{ all other $[,]$ trivial}\rangle.$$
These groups are nilpotent of class~$2$ and might be thought of as the direct product of~$n-1$
copies of the discrete Heisenberg group with one diagonal entry identified in
each copy. In Theorem~\ref{main theorem} we shall determine the subgroup zeta functions of
each~$G_n$, i.e. the Dirichlet series
$$\zeta_{G_n}^{\leq}(s):=\sum_{H\leq G_n} |G_n:H|^{-s}.$$
This generalizes the well-known formula (cf.~\cite{GSS/88}) for the
Heisenberg group~$G_2$
$$\zeta_{G_2}^{\leq}(s)=\zeta(s)\zeta(s-1)\zeta(2s-2)\zeta(2s-3)\zeta(3s-3)^{-1}$$
where $\zeta(s)=\sum_{i=1}^\infty i^{-s}$ denotes the Riemann zeta
function. 

The subgroup zeta function satisfies an Euler product
decomposition
$$\zeta_{G_n}^{\leq}(s)=\prod_{p \text{
    prime}}\zeta_{G_n,p}^{\leq}(s)$$
into local zeta functions
\begin{equation}
\zeta_{G_n,p}^{\leq}(s):=\sum_{H\leq_p G_n} |G_n:H|^{-s}\label{def
  local zeta function}
\end{equation}
enumerating subgroups of finite $p$-power index in~$G_n$.

\begin{theorem}\label{main theorem} For a prime~$p$
$$\zeta_{G_n,p}^{\leq}(s)=\zeta_p(s-1)\dots\zeta_p(s-n+1)\zeta_p(ns-n(n-1))W^\leq_{n-1}({\bf X},{\bf Y},p),$$
where 
\begin{equation}
W^\leq_{n-1}({\bf X},{\bf
  Y},p):=\sum_{I\subseteq\{1,\dots,n-2\}}b_{n-1,I}(p^{-1})W^\leq_{{n-1},I}({\bf
  X}, {\bf Y}), \label{subgroup sum}
\end{equation}
 with
 \begin{equation}
W^\leq_{{n-1},I}({\bf X}, {\bf Y}):=  \ul{X_0}\prod_{i\in I}\ul{X_i}+\ul{Y_{n-1}}\prod_{i\in I}\ul{Y_i}+\sum_{j\in
 I}\prod_{\substack{i\in I\\i\leq j}}\ul{Y_i}\prod_{\substack{i\in I\\i\geq j}}\ul{X_i}+\sum_{j\in I\cup\{n-1\}}\prod_{\substack{i\in I\\i< j}}\ul{Y_i}\prod_{\substack{i\in I\\i\geq j}}\ul{X_i},
\label{grenham formula}
\end{equation}
where $\ul{Z}:=\frac{Z}{1-Z}$ and the numerical data ${\bf X}=(X_0,\dots,X_{n-2})$, ${\bf
 Y}=(Y_1,\dots,Y_{n-1})$ is defined by
\begin{eqnarray}
Y_i &:=& p^{-(n-i)s+(n+i)(n-1-i)}, \quad 1\leq i \leq n-1, \label{num
  data Y} \\
X_i &:=& (p^{-2s+(n+1+i)})^{(n-1-i)}, \quad 0\leq i \leq n-2. \label{num
  data X}
\end{eqnarray}
The polynomial expression~$b_{n-1,I}(p)$ denotes the number of flags
of type~$I$ in $\Fp^{n-1}$ and $\zeta_p(s)=\sum_{i=0}^\infty p^{-is}$
denotes the $p$-th factor of the Riemann zeta function.
\end{theorem}

\begin{corollary}\label{coro funeq} For all $n\geq 2$ and all primes $p$ the
  following local functional equation holds:
\begin{equation}
\zeta_{G_n,p}^{\leq}(s)|_{p\rightarrow
  p^{-1}}=-p^{\binom{2n-1}{2}-(2n-1)s}\zeta_{G_n,p}^{\leq}(s).\label{grenham
  subgroup funeq}
\end{equation}
\end{corollary}

\bigskip
For completeness we record here the formulae for two other Dirichlet
series associated to the nilpotent groups~$G_n$, $n\geq2$. The zeta functions
\begin{eqnarray*}
\zeta^\triangleleft_{G_n}(s)&:=& \sum_{H\triangleleft G_n} |G_n:H|^{-s}\\
\zeta^\wedge_{G_n}(s)&:=& \sum_{\substack{H\leq G_n\\ \widehat{H}\isom \widehat{G_n}}} |G_n:H|^{-s}
\end{eqnarray*}
(here $\widehat{G}$ denotes the profinite completion of the group~$G$)
 both decompose as Euler products of {\sl local} zeta functions,
 defined in complete analogy to~(\ref{def local zeta function}). 

The explicit formulae and functional equations for
the factors of the {\sl normal} zeta functions~$\zeta^\triangleleft_{G_n}(s)$
were derived in~\cite{Voll/03a}:

\begin{theorem}{[\cite{Voll/03a}, Thm~5]}
 $$ \zeta^\triangleleft_{G_n,p}(s)=\prod_{i=0}^{n-1}\zeta_p(s-i)\; \zeta_p((2n-1)s-n(n-1))W^\triangleleft_{n-1}({\bf X},p),$$
where 
\begin{equation}
W^\triangleleft_{n-1}({\bf X},p):=
\sum_{I\subseteq\{1,\dots,n-2\}}b_{{n-1},I}(p^{-1})\prod_{i\in
  I}\frac{X_i}{1-X_i} \label{normal sum}
\end{equation}
and the numerical data ${\bf X}=(X_1,\dots,X_{n-2})$ is defined by
$$X_i=p^{-(2(n-i)-1)s+(n+i)(n-i-1)}\text{ for }i\in\{1,\dots,n-2\}.$$
For all $n\geq2$ and all primes~$p$ the following functional equation holds
\begin{equation}
\zeta^\triangleleft_{G_n,p}(s)|_{p\rightarrow
  p^{-1}}=-p^{\binom{2n-1}{2}-(3n-1)s}\zeta^\triangleleft_{G_n,p}(s).\label{grenham normal funeq} 
\end{equation}
\end{theorem}
 The functional equations~(\ref{grenham subgroup funeq})
 and~(\ref{grenham normal funeq}) were conjectured by du~Sautoy
 (Conjecture~$5.47$ in~\cite{duS-ennui/02}). The zeta functions~$\zeta^*_{G_n}(s)$,
$*\in\{\leq,\triangleleft\}$, $n=3,4,5$, were first computed by
D.~Grenham in his doctoral thesis~\cite{Grenham/88}. His computations are recorded
in~\cite{duSG-ghosts/00}.

The zeta functions~$\zeta^{\wedge}_{G_n}(s)$ were computed by
Mark Berman, a student of du~Sautoy, following the analysis of~\cite{duSLubotzky/96}.
\begin{proposition}{[Berman]} For all $n\geq 2$ and all primes~$p$
\begin{equation}
 \zeta^{\wedge}_{G_n,p}(s)=\zeta_p(ns-n(n-1))\prod_{i=0}^{n-2}\zeta_p(2s-n-1-i).\label{grenham
 hat formula}
\end{equation}
In particular
\begin{equation*}
\zeta^{\wedge}_{G_n,p}(s)|_{p\rightarrow
  p^{-1}}=(-1)^np^{5\binom{n}{2}-(3n-2)s}\zeta^{\wedge}_{G_n,p}(s).\label{grenham hat funeq} 
\end{equation*}
\end{proposition}

The Dirichlet series $\zeta^*_{G_n}(s)$, $*\in\{\leq,
\triangleleft,\wedge\}$ define analytic function on right
half-planes $\{s\in\mathbb{C}|\; \mathfrak{Re}(s)>\alpha_{G_n}^*\}$. Du
Sautoy and Grunewald showed in~\cite{duSG/00} that these abscissae of
convergence~$\alpha_{G_n}^*$ are always {\sl rational} numbers. But
computing them effectively is a difficult problem. So far
we had few examples for which they were not integers (the first one
being the abscissa of convergence~$\alpha^\leq_{F_{2,3}}$ of the
subgroup zeta function of~$F_{2,3}$, the free class-$2$-nilpotent
group on three generators, which was computed by
G. Taylor~\cite{Taylor/01}). And it is in fact not hard to see that 
$$\alpha^{\wedge}_{G_n}=\alpha^\triangleleft_{G_n}=n.$$
The existence of infinitely many non-integral rationals among the
abscissae of the subgroup zeta functions~$\zeta_{G_n,p}^{\leq}(s)$ is a
consequence of a result of Pirita Paajanen, a student of du~Sautoy,
obtained independently of Theorem~\ref{main theorem}:

\begin{proposition}{[Paajanen]}\label{proposition abscissa}
\begin{equation}
\alpha^\leq_{G_n}=\max_{1\leq l \leq n-2}\left\{
  n,\frac{(n+l)(n-l-1)+1}{n-l}\right\}.\label{abscissa formula}
\end{equation}
\end{proposition}
Note that the abscissa of convergence~$\alpha^\leq_{G_n}$ is not an integer
if~$2n-1$ is a prime~$\geq11$. Indeed, calculus shows that the right hand side
  of~(\ref{abscissa formula}) is greater than~$n$ for $n\geq 6$ and
  $\frac{(n+l)(n-l-1)+1}{n-l}\equiv \frac{1-2n}{n-l}\mod \Z$. (This
  observation was made by O. Sauzet.)

\bigskip
Interesting open questions concern the poles of the zeta
functions~$\zeta^*_{G_n}(s)$. It is not clear {\sl a priori} which factors
of the denominators of the rational functions~(\ref{subgroup sum}) and~(\ref{normal sum}) will cancel out when written in lowest
terms. Comparison of~(\ref{grenham
 hat formula}) with~(\ref{num  data X}) also yields that all the factors of the denominator
of~$\zeta^\wedge_{G_n}(s)$ occur as `candidate factors' in the
denominators of~(\ref{subgroup sum}). The example~$G_5$, however,
shows that they may fall victim to cancellation. 

\bigskip
The proofs of Theorem~\ref{main theorem} and its Corollary~\ref{coro
  funeq} will be given in Section~\ref{section proofs}. We will see
  that proving the functional equations~(\ref{grenham
  subgroup funeq}) (and~(\ref{grenham normal funeq})) reduces to establishing a combinatorial property of
  the rational functions~$W^\leq_{n-1}(\bfX,\bfY,p)$
  (and~$W^\triangleleft_{n-1}(\bfX,p)$) in~(\ref{grenham
  formula}) (and~(\ref{normal sum}), respectively) in Lemma~\ref{comb lemma 1}. To derive Theorem~\ref{main
  theorem} we analyse the Cartan decomposition of certain subquotients
  of lattices in the Lie ring associated to the nilpotent group~$G_n$.

\section{Proofs}\label{section proofs}
\subsection{Proof of Corollary~\ref{coro funeq}}
One checks easily that it is enough to show that
 
\begin{equation}
W^\leq_n({\bf X}^{-1},{\bf Y}^{-1},p^{-1})=(-1)^n p^{\binom{n}{2}}W^\leq_n({\bf X}, {\bf Y},p)\quad\forall n\geq1 \label{comb funeq},
\end{equation}
where of course $\bfX^{-1}=(X_0^{-1},\dots,X_{n-1}^{-1})$ etc.. (We
will omit the superscript $\leq$ from now on.) Note that the
functional equation~(\ref{comb funeq}) is independent of the actual numerical
data specified in~(\ref{num data X}) and~(\ref{num data Y}). To prove~(\ref{comb funeq})   we need two
combinatorial lemmata. Only Lemma~\ref{comb lemma 1}
depends on the actual shape of the rational functions~(\ref{grenham formula}).

We shall write $[{m}]$ for
$\{1,\dots,m\}$, $|I|$ for the cardinality of the finite
set~$I\subseteq[{m}]$, $I^c$ for $[{m}]\setminus I$
and $\uplus$ for the union of disjoint sets.

\begin{lemma}\label{comb lemma 1} 
$$
W_{n,J}(\bfX^{-1},\bfY^{-1})=(-1)^{|J|+1}\sum_{S\subseteq J}W_{n,S}(\bfX,\bfY)\quad\forall n\geq 1,\,J \subseteq [{n-1}].$$
\end{lemma}

\begin{proof} Note that $\ul{Z^{-1}}=-(1+\ul{Z})$. Thus
\begin{eqnarray*}
\lefteqn{(-1)^{|J|+1}W_{n,J}(\bfX^{-1},\bfY^{-1})=(1+\ul{X_0})\prod_{i\in
    J}(1+\ul{X_i})+(1+Y_n)\prod_{i\in J}(1+\ul{Y_i})}\\
&&\quad-\sum_{j\in J\cup\{n\}} \prod_{\substack{i \in
    J\\i< j}}(1+\ul{Y_i})\prod_{\substack{i\in J\\i\geq
    j}}(1+\ul{X_i})+\sum_{j\in J} (1+\ul{Y_j})\prod_{\substack{i \in J\\i<j}}(1+\ul{Y_i})\prod_{\substack{i\in J \\i\geq
    j}}(1+\ul{X_i})\\
&=&(1+\ul{X_0})\prod_{i\in
    J}(1+\ul{X_i})+(1+\ul{Y_n})\prod_{i\in J}(1+\ul{Y_i})-\prod_{i\in J}(1+\ul{Y_i})+\sum_{j \in J}\ul{Y_j}\prod_{\substack{i\in J\\ i<
    j}}(1+\ul{Y_i})\prod_{\substack{i\in J\\i\geq
    j}}(1+\ul{X_i})\\
&=&\sum_{S\subseteq J}\left(\ul{X_0}\prod_{s\in
    S}\ul{X_s}+\ul{Y_n}\prod_{s\in S}\ul{Y_s}+\sum_{j\in
    S}\prod_{\substack{s\in S\\s\leq j}}\ul{Y_s}\prod_{\substack{s\in
    S\\s\geq j}}\ul{X_s}+{\sum_{j\in S}\prod_{\substack{s\in S\\s\leq
    j}}\ul{Y_s}\prod_{\substack{s\in S\\s>j}}\ul{X_s}+\prod_{s\in
    S}\ul{X_s}}\right)\\
&=&\sum_{S\subseteq J}W_{n,S}(\bfX,\bfY),\text{ as desired.}
\end{eqnarray*}
\end{proof}
%_{=\sum_{j\in S\cup\{n\}}\prod_{\substack{s\in S\\s <   j}}\ul{Y_s}\prod_{\substack{s\in S\\s\geq j}}\ul{X_s}}
\begin{lemma} \label{comb lemma 2}
$$\sum_{J\supseteq I}W_{n,J}(\bfX^{-1},\bfY^{-1})=(-1)^n\sum_{J\supseteq
  I^c}W_{n,J}(\bfX,\bfY)\quad \forall n\geq 1,\; I\subseteq [{n-1}].$$
\end{lemma}
\begin{proof} By Lemma~\ref{comb lemma 1} we have
\begin{eqnarray*}
\sum_{J\supseteq I}W_{n,J}(\bfX^{-1},\bfY^{-1})&=&-\sum_{J\supseteq
  I}(-1)^{|J|}\sum_{S\subseteq J}W_{n,S}(\bfX,\bfY)\\
&=&-\sum_{R\subseteq[{n-1}]}c_RW_{n,R}(\bfX,\bfY),\text{ say,}
\end{eqnarray*}
where
\begin{eqnarray*}
c_R&=&\sum_{R\cup I\subseteq J}(-1)^{|J|}=(-1)^{|R\cup I|}\sum_{S\subseteq(R\cup I)^c}(-1)^{|S|}\\&=&\left.(-1)^{|R\cup I|}\cdot0^{|(R\cup I)^c|}=\left\{\begin{array}{cl}
(-1)^{n-1}&\mbox{if $R\supseteq I^c$,}\\
0&\mbox{otherwise.}
\end{array}\right.\right.
\end{eqnarray*}
\end{proof}

We are now ready to prove the functional equation~(\ref{comb funeq}). We will use the following notation:  Given an element~$w\in S_n$ of the
symmetric group on~$n$ letters $\{1,\dots,n\}$ we write 
$$\nu(w):=\{i\in[n-1]|\, w(i)>w(i+1)\}$$ for the
{\sl(descent) type} of~$w$. By $l(w)$ we denote the length of~$w$ as a
word in the standard generators for~$S_n$ and
by $w_0$ the longest element. We have the standard identities
\begin{eqnarray}
\nu(ww_0)&=& \nu(w)^c, \label{typeinversion}\\
l(w)+l(ww_0)&=&\binom{n}{2}. \label{length} 
\end{eqnarray}
For~$n\geq 1$ we have
\begin{eqnarray}
W_n({\bf X},{\bf
  Y},p)&=&\sum_{I\subseteq[{n-1}]}b_{n,I}(p^{-1})W_{n,I}({\bf X},{\bf
  Y})\nonumber\\
&=&\sum_{w\in S_n}p^{-l(w)}\sum_{J\supseteq \nu(w)}\wnjxy \label{5}.
\end{eqnarray}
We may rewrite the left hand side of~(\ref{comb
  funeq}) as
\begin{alignat*}{2}
W_n({\bf X}^{-1},{\bf Y}^{-1},p^{-1})&=p^{\binom{n}{2}}\sum_{w\in
  S_n}p^{-l(ww_0)}\sum_{J\supseteq \nu(w)}W_{n,J}(\bfX^{-1},\bfY^{-1})
 &\text{[by~(\ref{length}) and~(\ref{5})]}&\\
&=(-1)^np^{\binom{n}{2}}\sum_{w\in S_n}p^{-l(ww_0)}\sum_{J\supseteq
  \nu(ww_0)}W_{n,J}(\bfX,\bfY)&\text{
  [by~(\ref{typeinversion}) and Lemma~\ref{comb lemma 2}]}&\\
&=(-1)^np^{\binom{n}{2}}W_n({\bf X},{\bf Y},p), &\text{[by~(\ref{5})]}
\end{alignat*}
as desired.

\subsection{Proof of Theorem~\ref{main theorem}}
 We fix
  $n\geq2$ and omit the subscript~$n$. We write $Z=Z(G)$ for the centre
  of~$G$ and $\LieL:=G/Z\oplus Z$
for the Lie ring associated to~$G$ with Lie bracket induced from
taking commutators in~$G$. Let~$N$ denote the unique abelian ideal
  in~$\LieL$ with cyclic quotient (generated by~$x_nN$, say). We now fix a prime~$p$ and
denote by $\LieLp$, $N_p$, $Z_p$ the respective pro-$p$-completions.
Given an additive $\Zp$-submodule~$H\subseteq(\LieLp,+)$ we define the three
  subquotients of~$\LieLp$
\begin{eqnarray*}
H_1&:=&H\cap Z_p, \\
H_2&:=&(H\cap N_p)Z_p/Z_p,\\
H_3&:=&HN_p/ N_p \text{ with }|\LieLp:HN_p|=p^{m_3},\text{ say, }m_3\in \N\cup\{0\}.
\end{eqnarray*}
Lie bracketing with $x_n$ induces the linear isomorphism
\begin{eqnarray*}
\phi:N_p/Z_p &\rightarrow& Z_p \\
       x &\mapsto& (x,x_n).
\end{eqnarray*}
It is not hard to see that a lattice~$H\subseteq(\LieLp,+)$ is a subalgebra of~$\LieLp$ if
and only if $$p^{m_3}\phi(H_2)\subseteq H_1,\quad\text{ or, equivalently,
}\quad\phi(H_2)\subseteq  (p^{-m_3}H_1)\cap Z_p\;(\subseteq Z_p\otimes_{\Zp}\Qp).$$ We
say that $H_1\subseteq Z_p\isom \Zp^{n-1}$ is of {\sl
  type}~$\nu(H_1)=(I,{\bf r}_{I_0})$, where
$$I=\{i_1,\dots,i_l\}_<\subseteq[n-2]\footnote{ Here
  $\{i_1,\dots,i_l\}_<$ means $i_1<\dots< i_l$.},\;{\bf r}_{I_0}=(r_0,r_{i_1},\dots,r_{i_l})\in\N_{>0}\times\N_{\geq0}^{|I|}$$
if~$H_1$ has elementary divisor type\footnote{Note that
  in~\cite{Voll/03a} we could assume~$H_1$ to be maximal, i.e.~$r_0=0$.}
$$\left(\underbrace{p^{r_0},\dots,p^{r_0}}_{i_1},\underbrace{p^{r_0+r_{i_1}},\dots,p^{{r_0}+r_{i_1}}}_{i_2-i_1},\dots,\underbrace{p^{\sum_{i\in  I_0}r_i},\dots,p^{\sum_{i\in      I_0}r_i}}_{n-1-i_{l}}\right) $$
where we wrote $I_0$ for $I\cup {\{0\}}$. One checks easily that
\begin{eqnarray*}
\zeta^\leq_{G,p}(s)&=&\sum_{\substack{(H_1,H_2,H_3)\\ \phi(H_2)\subseteq
    (p^{-m_3}H_1)\cap Z_p}}|Z_p:H_1|^{n-s}||N_p/Z_p:H_2|^{1-s}|\LieLp/N_p:H_3|^{-s}\\
&=&\prod_{i=1}^{n-1}\zeta_p(s-i)\cdot\sum_{\substack{(H_1,H_2,H_3)\\\phi(H_2)=   (p^{-m_3}H_1)\cap Z_p}}|Z_p:H_1|^{n-s}||N_p/Z_p:H_2|^{1-s}p^{-m_3s}\\
&=&\prod_{i=1}^{n-1}\zeta_p(s-i)\cdot\underbrace{\sum_{H_1\subseteq Z_p}|Z_p:H_1|^{n-s}\underbrace{\sum_{m_3\geq
    0}|Z_p:(p^{-m_3}H_1)\cap
    Z_p|^{1-s}p^{-m_3s}}_{=:Z(I,{\bf r}_{I_0})}}_{=:Z(s)}
\end{eqnarray*}
where $H_1,H_2$, $H_3$ run through the sublattices of
$Z_p,N_p/Z_p$ and $\LieLp/N_p$, respectively. We observe that $Z(I,{\bf r}_{I_0})$ is well-defined as the
respective sum only depends on the type of~$H_1$. Recall
(e.g. from~\cite{Voll/03a}, Lemma~1) that 
\begin{equation*}f(I,{\bf r}_{I_0 },p):=|\{ H_1\subseteq
  Z_p|\;\nu(H_1)=(I,{\bf r}_{{ I_0}})\}|=b_{n-1,I}(p^{-1})p^{\sum_{i\in I_0}r_i(n-1-i)i},
\end{equation*} 
and note that
\begin{equation*}|Z_p:H_1|=p^{\sum_{i\in
      I_0}r_i(n-1-i)}.
\end{equation*}
Therefore
\begin{eqnarray}
Z(s)&=&\sum_{I \subseteq[n-2]}b_{n-1,I}(p^{-1})\sum_{{\bf r}_{I_0}}p^{\sum_{i\in I_0}r_i(n-1-i)i}p^{(n-s)\sum_{i\in
    I_0}r_i(n-1-i)} Z(I,{\bf r}_{I_0})\nonumber\\
&=&\sum_{I\subseteq[n-2]}b_{n-1,I}(p^{-1})\sum_{{\bf r}_{I_0}}\prod_{j\in I_0}y_j^{r_j}Z(I,{\bf r}_{I_0})\label{Z(s)}
\end{eqnarray}
where we set
$$y_i:=p^{(n-s+i)(n-1-i)}  \text{ for }i \in \{0,1,\dots,n-2\}. $$ 
To compute the sum $Z(I,{\bf r}_{I_0})$ it will prove advantageous to decompose its range of summation as
$$\{m_3\in\N_{\geq0}\}=\biguplus_{j\in I_0}\mathcal{M}_j({\bf r}_{I_0})\uplus\{m_3\geq
 \sum_{i\in I_0}r_i\},\text{ say,}$$
where $$\mathcal{M}_j({\bf r}_{I_0}):= \left\{{\sum_{{i\in I_0,\;i<j}}r_i\leq m_3<
  \sum_{{i\in I_0,\;i\leq j}}r_i}\right\}.$$
Notice that only $\mathcal{M}_0$ may be empty and that
\begin{equation*}|Z_p:(p^{-m_3}H_1)\cap Z_p|=1\text{ if and only if }m_3\geq
 \sum_{i\in I_0}r_i.
\end{equation*} We may now write 
 \begin{eqnarray*}
Z(I,{\bf r}_{I_0})&=&\sum_{j\in I_0}\sum_{m_3\in \mathcal{M}_j({\bf r}_{I_0})}|Z_p:(p^{-m_3}H_1)\cap
    Z_p|^{1-s}p^{-m_3s}+\sum_{m_3\geq \sum_{i\in
  I_0}r_i}p^{-m_3s}\\
&=&\sum_{j\in I_0}p^{(-s)\sum_{\substack{i\in
  I_0\\i<j}}r_i+(1-s)\sum_{\substack{i\in I_0\\i\geq
  j}}r_i(n-1-i)}\left(\frac{1-\left(\frac{p^{-s}}{p^{(1-s)(n-1-j)}}\right)^{r_j}}{1-\left(\frac{p^{-s}}{p^{(1-s)(n-1-j)}}\right)}\right)
 +\sum_{m_3\geq \sum_{i\in I_0}r_i}p^{-m_3s}\\
&=&\sum_{j\in  I_0} \prod_{\substack{i\in I_0\\i<j}}z^{r_i}\prod_{\substack{i\in
  I_0\\i\geq
  j}}x_i^{r_i}\left(\frac{1-\left(\frac{z}{x_j}\right)^{r_j}}{1-\frac{z}{x_j}}\right)+\frac{z^{
  \sum_{i\in I_0}r_i}}{1-z}
\end{eqnarray*}
with the following abbreviations
\begin{eqnarray*}
x_i&:=&p^{(1-s)(n-1-i)} \text{ for }i \in \{0,1,\dots,n-2\}\\
%y_i&:=&p^{(n-s+i)(n-1-i)}  \text{ for }i \in \{0,1,\dots,n-1\}\\
z &:=& p^{-s}.
\end{eqnarray*}
We rewrite~(\ref{Z(s)}) as
\begin{eqnarray}
\lefteqn{Z(s)=\frac{1}{1-y_0z}\cdot}\nonumber\\
&&\sum_{I\subseteq[n-2]}b_{n-1,I}(p^{-1})\left(\prod_{i\in
  I_0}\frac{x_iy_i}{1-x_iy_i}+\frac{1}{1-y_jz}\prod_{\substack{i\in I_0\\i<j}}\frac{y_iz
}{1-y_iz}\prod_{\substack{i\in
  I_0\\i\geq
  j}}\frac{x_iy_i}{1-x_iy_i}+\frac{1}{1-z}\prod_{i\in I}\frac{y_iz
}{1-y_iz}\right)\label{almost final}
\end{eqnarray}
as one checks without difficulty. The result now follows by setting
\begin{eqnarray*}
Y_i&:=&y_iz,\\
X_i&:=&x_iy_i\text{ for }i \in \{0,1,\dots,n-2\}
\end{eqnarray*}
and rearranging the sum in brackets\footnote{Paajanen's result (Proposition~\ref{proposition abscissa}) amounts to
showing that replacing the sum in brackets in~(\ref{almost final})
by its last summand yields a Dirichlet series whose absissa of
convergence equals~$\alpha^\leq_{G_n}$.} in~(\ref{almost final}).

\begin{acknowledgements} We should like to thank the
  UK's Engineering and Physical Sciences Research Council~(EPSRC) for
  their support in the form of a Postdoctoral Fellowship. 
\end{acknowledgements}

\bibliographystyle{amsplain}
\bibliography{thebibliography}

\end{document}